\documentstyle[12pt]{article}

\begin{document}

\title{{\bf Periodical solutions of Poisson-gradient dynamical systems with
periodical potential}}
\author{Constantin Udri\c ste, Iulian Duca}
\date{}
\maketitle

\begin{abstract}
The main purpose of this paper is the study of the action that produces
Poisson-gradient systems and their multiple periodical solutions. The
Section 1 establishes the basic tools.

The section 2 underlines conditions in which the action $\varphi \left(
u\right) = \displaystyle\displaystyle\int_{T_{0}}\left[ \displaystyle%
\displaystyle\frac{1}{2}\left| \displaystyle\displaystyle\frac{\partial u}{%
\partial t}\right| ^{2}+F\left( t,u\left( t\right) \right) \right]
dt^{1}\wedge ...\wedge dt^{p}$, that produces the Poisson-gradient systems,
is continuous, and some conditions in which the general action $\varphi
\left( u\right) = \displaystyle\displaystyle\int_{T_{0}}L\left( t,u\left(
t\right), \displaystyle\displaystyle\frac{\partial u}{\partial t}\left(
t\right) \right) dt^{1}\wedge ...\wedge dt^{p}$ is continuously
differentiable.

The Section 3 studies the multiple periodical solutions of a
Poisson-gradient system in the case when the potential function $F$ has a
spatial periodicity.
\end{abstract}



\section{Introduction}

\hspace{0cm} In this paper we will note by $W_{T}^{1,2}$ the Sobolev spaces
of the $u\in L^{2}\left[ T_{0},R^{n}\right] $ functions, which have the weak
derivative $\displaystyle\frac{\partial u}{\partial t}\in L^{2}\left[
T_{0},R^{n}\right] $, $T_{0}=\left[ 0,T^{1}\right] \times ...\times \left[
0,T^{p}\right] \subset R^{p}.$ The weak derivatives are defined using the
space $C_{T}^{\infty }$\ of all indefinitely differentiable multiple $T$%
-periodic function from $R^{p}$ into $R^{n}$. We consider $H_{T}^{1}$ \ the
Hilbert space of the $W_{T}^{1,2.}$. The geometry on $H_{T}^{1}$ is realized
by the scalar product 
\[
\left\langle u,v\right\rangle =\displaystyle\int_{T_{0}}\left( \delta
_{ij}u^{i}\left( t\right) v^{j}\left( t\right) +\delta _{ij}\delta ^{\alpha
\beta }\displaystyle\frac{\partial u^{i}}{\partial t^{\alpha }}\left(
t\right) \displaystyle\frac{\partial v^{j}}{\partial t^{\beta }}\left(
t\right) \right) dt^{1}\wedge ...\wedge dt^{p}, 
\]
and the associated Euclidean norm $\left\| .\right\| .$ These are induced by
the scalar product (Riemannian metric) 
\[
G=\left( 
\begin{array}{cc}
\delta _{ij} & 0 \\ 
0 & \delta ^{\alpha \beta }\delta _{ij}
\end{array}
\right) 
\]
on $R^{n+np}$ (see the jet space $J^{1}\left( T_{0},R^{n}\right) $).

Let $t=\left( t^{1},...,t^{p}\right) $ be a generic point in $R^{p}$. Then
the opposite faces of the parallelepiped $T_{0}$ can be described by the
equations 
\[
S_{i}^{-}:t^{i}=0,S_{i}^{+}:t^{i}=T^{i} 
\]
for each $i=1,...,p.$ We denote 
\[
\left\| u\right\| _{L^{2}}=\displaystyle\int_{T_{0}}\delta _{ij}u^{i}\left(
t\right) v^{j}\left( t\right) dt^{1}\wedge ...\wedge dt^{p}, 
\]
\[
\left\| \frac{\partial u}{\partial t}\right\| _{L^{2}}=\displaystyle%
\int_{T_{0}}\delta _{ij}\delta ^{\alpha \beta }\displaystyle\frac{\partial
u^{i}}{\partial t^{\alpha }}\left( t\right) \displaystyle\frac{\partial v^{j}%
}{\partial t^{\beta }}\left( t\right) dt^{1}\wedge ...\wedge dt^{p}, 
\]
\[
\left( u,v\right) =\delta _{ij}u^{i}v^{j},\left| u\right| =\sqrt{\delta
_{ij}u^{i}u^{j}}. 
\]

We study the extremals of the action 
\[
\varphi \left( u\right) =\displaystyle\int_{T_{0}}\left[ \displaystyle\frac{1%
}{2}\left| \displaystyle\frac{\partial u}{\partial t}\right| ^{2}+F\left(
t,u\left( t\right) \right) \right] dt^{1}\wedge ...\wedge dt^{p} 
\]
on $H_{T}^{1}$ in the case when the potential $F$ function has spatial
periodicity. So, we consider that there exist $P_{1},....P_{n}\in R$ so that 
$F\left( t,x+P_{i}e_{i}\right) =F\left( t,x\right) $, for any $t\in T_{0}$, $%
x\in R^{n}$ and any $i\in \{1,...n\}$. The vectors $e_{1},...,e_{n}$ create
a canonical base\ in the Euclidian space $R^{n}.$ The extremals of the
action $\varphi $ are being determined with the minimizing sequences method.
For the existence of the minimizing bounded sequence we need to introduce
the spatial periodicity condition of the potential function $F$. The
function that realizes the minimum of the action $\varphi $ verifies the
Euler-Lagrange equations, which in the case of the Lagrangian 
\[
L\left( t,u\left( t\right) ,\displaystyle\frac{\partial u}{\partial t}%
\right) =\displaystyle\frac{1}{2}\left| \displaystyle\frac{\partial u}{%
\partial t}\right| ^{2}+F\left( t,u\left( t\right) \right) 
\]
defined on $H_{T}^{1}$, reduces to a Poisson-gradient $PDE_{S}$, $\Delta
u\left( t\right) =\nabla F\left( t,u\left( t\right) \right) ,$ 
\[
u\mid _{S_{i}^{-}}=u\mid _{S_{i}^{+}},\frac{\partial u}{\partial t}\mid
_{S_{i}^{-}}=\frac{\partial u}{\partial t}\mid _{S_{i}^{+}},i=1,...,p. 
\]

\section{The action that produces Poisson-gradient \newline
systems}

\subsection{Multi-time Euler-Lagrange equations}

We consider the multi-time variable $t=\left( t^{1},...,t^{p}\right) \in
R^{p}$, the functions $x^{i}:R^{p}\rightarrow R,\left(
t^{1},...,t^{p}\right) \rightarrow $ $x^{i}\left( t^{1},...,t^{p}\right) $, $%
i=1,...n,$ and we denote $x_{\alpha }^{i}=\displaystyle\frac{\partial x^{i}}{%
\partial t^{\alpha }},\alpha =1,...,p$. The Lagrange function 
\[
L:R^{p+n+np}\rightarrow R,\left( t^{\alpha },x^{i},x_{\alpha }^{i}\right)
\rightarrow L\left( t^{\alpha },x^{i},x_{\alpha }^{i}\right) 
\]
gives the Euler-Lagrange equations 
\[
\displaystyle\frac{\partial }{\partial t^{\alpha }}\displaystyle\frac{%
\partial L}{\partial x_{\alpha }^{i}}=\displaystyle\frac{\partial L}{%
\partial x^{i}},\;i=1,...,n,\;\alpha =1,...p 
\]
(second order $PDEs$ system on the n-dimensional space).

The multi-time Lagrangian and Hamiltonian dynamics is based on the concept
of multisymplecticity (polysymplecticity) [3], [5]-[12]. Our task is to
develope some ideas from [1]-[2], [12]-[14] having in mind the single-time
theory in [4].

\subsection{Continuous action}

We consider the Lagrange function $L:T_{0}\times R^{n}\times
R^{np}\rightarrow R,\left( t^{\alpha },u^{i},u_{\alpha }^{i}\right)
\rightarrow L\left( t^{\alpha },u^{i},u_{\alpha }^{i}\right) ,$ $u_{\alpha
}^{i}=\displaystyle\frac{\partial u^{i}}{\partial t^{\alpha }},\alpha
=1,...,p,i=1,...,n,$ $u^{i}:T_{0}\rightarrow R,\left( t^{1},...,t^{p}\right)
\rightarrow u^{i}\left( t^{1},...,t^{p}\right) ,$ 
\[
L\left(t^{\alpha },u^{i},u_{\alpha }^{i}\right) =\displaystyle\frac{1}{2}%
\left| \displaystyle\frac{\partial u}{\partial t}\right| ^{2}+F\left(
t,u\left( t\right) \right) . 
\]

In the following result we will establish the conditions in which the action 
\[
\varphi \left( u\right) =\displaystyle\int_{T_{0}}L\left( t,u\left( t\right)
,\displaystyle\frac{\partial u}{\partial t}\left( t\right) \right)
dt^{1}\wedge ...\wedge dt^{p} 
\]
is continuous.

{\bf Theorem 1}. {\it Let $F:T_{0}\times R^{n}\rightarrow R$}${\it ,}${\it \ 
$\left( t,u\right) \rightarrow F\left( t,u\right) $\ be a measurable
function in $t$ for any $u\in R^{n}$ and continuously differentiable in $u$
for any $t$ $\in T_{0}$}${\it ,}${\it \ $T_{0}=\left[ 0,T^{1}\right] \times
...\times \left[ 0,T^{p}\right] \subset R^{p}$}${\it .}${\it \ If exists $%
M\geq 0$ and $g\in C\left( T_{0},R\right) $ such that $\left| \nabla
_{u}F\left( t,u\right) \right| \leq M\left| u\right| +g\left( t\right) $ for
any $t$ $\in T_{0}$ and any $u\in R$}${\it ,}${\it \ then $\varphi \left(
u\right) =\displaystyle\int_{T_{0}}\left[ \displaystyle\frac{1}{2}\left| %
\displaystyle\frac{\partial u}{\partial t}\right| ^{2}+F\left( t,u\left(
t\right) \right) \right] dt^{1}\wedge ...\wedge dt^{p}$ is continuous in $%
H_{T}^{1}$}${\it .}$

{\bf Proof}. We consider the $\left( u_{k}\right) _{k\in N}$ sequence which
is convergent in $H_{T}^{1}$ and we will note by $u$ his limit. This leads
us to the fact that 
\[
\left\| u_{k}-u\right\| ^{2}=\displaystyle\int_{T_{0}}\left[ \left|
u_{k}\left( t\right) -u\left( t\right) \right| ^{2}+\left| \displaystyle%
\frac{\partial u_{k}}{\partial t}\left( t\right) -\displaystyle\frac{%
\partial u}{\partial t}\left( t\right) \right| ^{2}\right] dt^{1}\wedge
...\wedge dt^{p} 
\]
\[
=\left\| u_{k}-u\right\| _{L^{2}}^{2}+\left\| \displaystyle\frac{\partial
u_{k}}{\partial t}-\displaystyle\frac{\partial u}{\partial t}\right\|
_{L^{2}}^{2}\rightarrow 0 
\]
when $k\rightarrow \infty $. The convergence of $u_{k}$ to $u$ in $H_{T}^{1}$
is equivalent to the convergence of $u_{k}$ to $u$ and to the convergence of 
$\displaystyle\frac{\partial u_{k}}{\partial t}$ to $\displaystyle\frac{%
\partial u}{\partial t}$ in $L^{2}$. By consequence $\left\| u_{k}\right\|
_{L^{2}}$ and $\left\| \displaystyle\frac{\partial u_{k}}{\partial t}%
\right\| _{L^{2}}$ are bounded in $R$. In order to show the continuity of $%
\varphi $, we will prove that $\left| \varphi \left( u_{k}\right) -\varphi
\left( u\right) \right| \rightarrow 0$ when $u_{k}\rightarrow u$ in $%
H_{T}^{1}$. In the evaluations that we will do for $\left| \varphi \left(
u_{k}\right) -\varphi \left( u\right) \right| $ we will utilize the \
following inequality $\left| \left\| u\right\| ^{2}-\left\| v\right\|
^{2}\right| \leq \left\| u\right\| \left\| u-v\right\| +\left\| v\right\|
\left\| u-v\right\| $ which is true in any vectorial space with an scalar
product. So, we obtain: 
\[
\begin{array}{l}
\left| \varphi \left( u_{k}\right) -\varphi \left( u\right) \right| \leq %
\displaystyle\frac{1}{2}\left| \displaystyle\int_{T_{0}}\left| \displaystyle%
\frac{\partial u_{k}}{\partial t}\right| ^{2}dt^{1}\wedge ...\wedge dt^{p}-%
\displaystyle\int_{T_{0}}\left| \displaystyle\frac{\partial u}{\partial t}%
\right| ^{2}dt^{1}\wedge ...\wedge dt^{p}\right| \\ 
\noalign{\medskip}+\displaystyle\int_{T_{0}}\left| F\left( t,u_{k}\left(
t\right) \right) -F\left( t,u\left( t\right) \right) \right| dt^{1}\wedge
...\wedge dt^{p}=\displaystyle\frac{1}{2}\left| \left\| \displaystyle\frac{%
\partial u_{k}}{\partial t}\right\| _{L^{2}}^{2}-\left\| \displaystyle\frac{%
\partial u}{\partial t}\right\| _{L^{2}}^{2}\right| \\ 
\displaystyle+\int_{T_{0}}\left| \displaystyle\int_{0}^{1}\left( \nabla
_{u}F\left( t,u_{k}\left( t\right) +s\left( u\left( t\right) -u_{k}\left(
t\right) \right) \right) ,u\left( t\right) -u_{k}\left( t\right) \right)
ds\right| dt^{1}\wedge ...\wedge dt^{p}
\end{array}
\]
\[
\begin{array}{l}
\displaystyle\leq \frac{1}{2}\left[ \left\| \displaystyle\frac{\partial u_{k}%
}{\partial t}\right\| _{L^{2}}\left\| \displaystyle\frac{\partial u_{k}}{%
\partial t}-\displaystyle\frac{\partial u}{\partial t}\right\|
_{L^{2}}+\left\| \displaystyle\frac{\partial u}{\partial t}\right\|
_{L^{2}}\left\| \displaystyle\frac{\partial u_{k}}{\partial t}-\displaystyle%
\frac{\partial u}{\partial t}\right\| _{L^{2}}\right] \\ 
\displaystyle+\int_{T_{0}}\left( \displaystyle\int_{0}^{1}\left| \nabla
_{u}F\left( t,u_{k}\left( t\right) +s\left( u\left( t\right) -u_{k}\left(
t\right) \right) \right) \right| \left| u\left( t\right) -u_{k}\left(
t\right) \right| ds\right) \\ 
dt^{1}\wedge ...\wedge dt^{p} \\ 
\displaystyle\leq \frac{1}{2}\left[ \left\| \displaystyle\frac{\partial u_{k}%
}{\partial t}\right\| _{L^{2}}\left\| \displaystyle\frac{\partial u_{k}}{%
\partial t}-\displaystyle\frac{\partial u}{\partial t}\right\|
_{L^{2}}+\left\| \displaystyle\frac{\partial u}{\partial t}\right\|
_{L^{2}}\left\| \displaystyle\frac{\partial u_{k}}{\partial t}-\displaystyle%
\frac{\partial u}{\partial t}\right\| _{L^{2}}\right] \\ 
\displaystyle+\int_{T_{0}}\left( \displaystyle\int_{0}^{1}\left( M\left|
u_{k}\left( t\right) +s\left( u\left( t\right) -u_{k}\left( t\right) \right)
\right| +g\left( t\right) \right) \left| u\left( t\right) -u_{k}\left(
t\right) \right| ds\right) \\ 
dt^{1}\wedge ...\wedge dt^{p} \\ 
\displaystyle\leq \frac{1}{2}\left[ \left\| \displaystyle\frac{\partial u_{k}%
}{\partial t}\right\| _{L^{2}}\left\| \displaystyle\frac{\partial u_{k}}{%
\partial t}-\displaystyle\frac{\partial u}{\partial t}\right\|
_{L^{2}}+\left\| \displaystyle\frac{\partial u}{\partial t}\right\|
_{L^{2}}\left\| \displaystyle\frac{\partial u_{k}}{\partial t}-\displaystyle%
\frac{\partial u}{\partial t}\right\| _{L^{2}}\right] \\ 
\displaystyle+\int_{T_{0}}\left( \displaystyle\int_{0}^{1}\left( M\left|
su\left( t\right) +\left( 1-s\right) \left( u_{k}\left( t\right) -u\left(
t\right) \right) \right| +g\left( t\right) \right) u\left( t\right) \right.
\left. \left| -u_{k}\left( t\right) \right| ds\right) \\ 
dt^{1}\wedge ...\wedge dt^{p}.
\end{array}
\]
Because the sequence $\left\| \displaystyle\frac{\partial u_{k}}{\partial t}%
\right\| _{L^{2}}$ is bounded, it exists $C_{1}$ such that $\left\| %
\displaystyle\frac{\partial u_{k}}{\partial t}\right\| _{L^{2}}\leq C_{1}$ \
for any $k\in N$. In the following, we will note by $C_{2}=\displaystyle%
\max_{t\in T_{0}}g\left( t\right) $ and we have 
\[
\begin{array}{l}
\left| \varphi \left( u_{k}\right) -\varphi \left( u\right) \right| \leq %
\displaystyle\frac{1}{2}C_{1}\left\| \displaystyle\frac{\partial u_{k}}{%
\partial t}-\displaystyle\frac{\partial u}{\partial t}\right\| _{L^{2}}+%
\displaystyle\frac{1}{2}\left\| \displaystyle\frac{\partial u}{\partial t}%
\right\| _{L^{2}}\left\| \displaystyle\frac{\partial u_{k}}{\partial t}-%
\displaystyle\frac{\partial u}{\partial t}\right\| _{L^{2}} \\ 
\noalign{\medskip}+M\left\| u\right\| _{L^{2}}\left( \displaystyle%
\int_{T_{0}}\left| u\left( t\right) -u_{k}\left( t\right) \right|
^{2}dt^{1}\wedge ...\wedge dt^{p}\right) ^{\displaystyle\frac{1}{2}}
\end{array}
\]
\[
\begin{array}{l}
\noalign{\medskip}+M\displaystyle\int_{T_{0}}\left| u_{k}\left( t\right)
-u\left( t\right) \right| ^{2}dt^{1}\wedge ...\wedge dt^{p}+C_{2}%
\displaystyle\int_{T_{0}}\left| u\left( t\right) -u_{k}\left( t\right)
\right| dt^{1}\wedge ...\wedge dt^{p} \\ 
\noalign{\medskip}\displaystyle\leq \frac{1}{2}\left( C_{1}+\left\| %
\displaystyle\frac{\partial u}{\partial t}\right\| _{L^{2}}\right) \left\| %
\displaystyle\frac{\partial u_{k}}{\partial t}-\displaystyle\frac{\partial u%
}{\partial t}\right\| _{L^{2}} \\ 
\noalign{\medskip}+\left( M\left\| u\right\| _{L^{2}}+C_{2}\left(
T^{1}...T^{p}\right) ^{\displaystyle\frac{1}{2}}\right) \ \left( %
\displaystyle\int_{T_{0}}\left| u\left( t\right) -u_{k}\left( t\right)
\right| ^{2}dt^{1}\wedge ...\wedge dt^{p}\right) ^{\displaystyle\frac{1}{2}}
\end{array}
\]

\[
\begin{array}{l}
\noalign{\medskip}+M\displaystyle\int_{T_{0}}\left| u_{k}\left( t\right)
-u\left( t\right) \right| ^{2}dt^{1}\wedge ...\wedge dt^{p}=\displaystyle%
\frac{1}{2}\left( C_{1}+\left\| \displaystyle\frac{\partial u}{\partial t}%
\right\| _{L^{2}}\right) \left\| \displaystyle\frac{\partial u_{k}}{\partial
t}-\displaystyle\frac{\partial u}{\partial t}\right\| _{L^{2}} \\ 
\noalign{\medskip}+\left( M\left\| u\right\| _{L^{2}}+C_{2}\left( T^{1}\cdot
\cdot \cdot T^{p}\right) ^{\displaystyle\frac{1}{2}}\right) \left\|
u-u_{k}\right\| _{L^{2}}+M\left\| u_{k}-u\right\| _{L^{2}}^{2}.
\end{array}
\]
Because $\left\| \displaystyle\frac{\partial u_{k}}{\partial t}-\displaystyle%
\frac{\partial u}{\partial t}\right\| _{L^{2}}\rightarrow 0$ and $\left\|
u_{k}-u\right\| _{L^{2}}\rightarrow 0$, when $k\rightarrow \infty $ it
results that $\varphi \left( u_{k}\right) \rightarrow \varphi \left(
u\right) $; from here we obtain the continuity of $\varphi $ in $H_{T}^{1}$.

\subsection{Continuously differentiable action}

In order to obtain a more general result then the one found in the previous
theorem , we define the action 
\[
\varphi :W_{T}^{1,2}\rightarrow R,\varphi \left( u\right) =\displaystyle%
\int_{T_{0}}L\left( t,u\left( t\right) ,\displaystyle\frac{\partial u}{%
\partial t}\left( t\right) \right) dt^{1}\wedge ...\wedge dt^{p},\; 
\]
\[
T_{0}=\left[ 0,T^{1}\right] \times ...\times \left[ 0,T^{p}\right] \subset
R^{p}. 
\]
Concerning this action we have the following Theorem which extends the
particular case $p=1$ from [4].

\bigskip {\bf Theorem 2}. {\it We consider $L:T_{0}\times R^{n}\times
R^{np}\rightarrow R,\left( t,x,y\right) \rightarrow L\left( t,x,y\right) $,
a measurable function in $t$ for any $\left( x,y\right) \in R^{n}\times
R^{np}$ and with the continuous partial derivatives in x and y for any $t\in
T_{0}$. If here exist $a\in C^{1}\left( R^{+},R^{+}\right) $ with the
derivative $a^{\prime }$\ bounded from above, $b\in C\left(
T_{0},R^{n}\right) $ such that for any $t\in T_{0}$ and any $\left(
x,y\right) \in R^{n}\times R^{np}$ to have } 
$$
\begin{array}{l}
\left| L\left( t,x,y\right) \right| \leq a\left( \left| x\right| +\left|
y\right| ^{2}\right) b\left( t\right) , \\ 
\noalign{\medskip}\left| \nabla _{x}L\left( t,x,y\right) \right| \leq
a\left( \left| x\right| \right) b\left( t\right) , \\ 
\noalign{\medskip}\left| \nabla _{y}L\left( t,x,y\right) \right| \leq
a\left( \left| y\right| \right) b\left( t\right) ,
\end{array}
\eqno(1) 
$$
{\it then, the functional $\varphi $ has continuous partial derivatives in $%
W_{T}^{1,2}$ and his gradient derives from the formula } 
$$
\begin{array}{lcl}
{\it \left( \nabla \varphi \left( u\right) ,v\right) } & {\it =} & {\it %
\displaystyle\int_{T_{0}}\left[ \left( \nabla _{x}L\left( t,u\left( t\right)
,\displaystyle\frac{\partial u}{\partial t}\right) ,v\left( t\right) \right)
\right. } \\ 
\noalign{\medskip} & {\it +} & {\it \left. \left( \nabla _{y}L\left(
t,u\left( t\right) ,\displaystyle\frac{\partial u}{\partial t}\left(
t\right) \right) ,\displaystyle\frac{\partial v}{\partial t}\left( t\right)
\right) \right] dt^{1}\wedge ...\wedge dt^{p}.}
\end{array}
\eqno(2) 
$$

{\bf Proof}. It is enough to prove that $\varphi $ has the derivative $%
\varphi ^{\prime }\left( u\right) \in \left( W_{T}^{1,2}\right) ^{\ast }$
given by the relation $\left( 2\right) $ and the function $\varphi ^{\prime
}:W_{T}^{1,2}\rightarrow \left( W_{T}^{1,2}\right) ^{\ast }$, $u\rightarrow
\varphi ^{\prime }\left( u\right) $ is continuous. We consider $u,v\in
W_{T}^{1,2}$, $t\in T_{0}$, $\lambda \in \left[ -1,1\right] $. We build the
functions 
\[
F\left( \lambda ,t\right) =L\left( t,u\left( t\right) +\lambda v\left(
t\right) ,\frac{\partial u}{\partial t}\left( t\right) +\lambda \frac{%
\partial v}{\partial t}\left( t\right) \right) 
\]
and 
\[
\Psi \left( \lambda \right) =\int_{T_{0}}F\left( \lambda ,t\right)
dt^{1}\wedge ...\wedge dt^{p}. 
\]
Because the derivative $a^{\prime }$ is bounded from above, exist\ $M>0$
such that $\displaystyle\frac{a\left( \left| u\right| \right) -a\left(
0\right) }{\left| u\right| }=a^{\prime }\left( c\right) \leq M.$ This means
that $a\left( \left| u\right| \right) \leq M\left| u\right| +a\left(
0\right) .$ On the other side

\[
\frac{\partial F}{\partial \lambda }\left( \lambda ,t\right) =\left( \nabla
_{x}L\left( t,u\left( t\right) +\lambda v\left( t\right) ,\frac{\partial u}{%
\partial t}\left( t\right) +\lambda \frac{\partial v}{\partial t}\left(
t\right) \right) ,v\left( t\right) \right) 
\]
\[
+\left( \nabla _{y}L\left( t,u\left( t\right) +\lambda v\left( t\right) ,%
\frac{\partial u}{\partial t}\left( t\right) +\lambda \frac{\partial v}{%
\partial t}\left( t\right) \right) ,\frac{\partial v}{\partial t}\left(
t\right) \right) \leq a\left( \left| u\left( t\right) +\lambda v\left(
t\right) \right| \right) 
\]
\[
b\left( t\right) \left| v\left( t\right) \right| +a\left( \left| \frac{%
\partial u}{\partial t}\left( t\right) +\lambda \frac{\partial v}{\partial t}%
\left( t\right) \right| \right) b\left( t\right) \left| \frac{\partial v}{%
\partial t}\left( t\right) \right| 
\]
\[
\leq b_{0}\left( M\left( \left| u\left( t\right) \right| +\left| v\left(
t\right) \right| \right) +a\left( 0\right) \right) \left| v\left( t\right)
\right| 
\]
\[
+b_{0}\left( M\left( \left| \frac{\partial u}{\partial t}\left( t\right)
\right| +\left| \frac{\partial v}{\partial t}\left( t\right) \right| \right)
+a\left( 0\right) \right) \left| \frac{\partial v}{\partial t}\left(
t\right) \right| , 
\]
where 
\[
b_{0}=\displaystyle\max_{t\in T_{0}}b\left( t\right) . 
\]
Then, we have $\left| \displaystyle\frac{\partial F}{\partial \lambda }%
\left( \lambda ,t\right) \right| \leq d\left( t\right) \in L^{1}\left(
T_{0},R^{+}\right) $. Then Leibniz formula of differentiation under integral
sign is applicable and 
\[
\frac{\partial \Psi }{\partial \lambda }\left( 0\right) =\int_{T_{0}}\frac{%
\partial F}{\partial \lambda }\left( 0,t\right) dt^{1}\wedge ...\wedge
dt^{p}=\int_{T_{0}}\left[ \left( \nabla _{x}L\left( t,u\left( t\right) ,%
\frac{\partial u}{\partial t}\left( t\right) \right) ,v\left( t\right)
\right) \right. 
\]
\[
\left. +\left( \nabla _{y}L\left( t,u\left( t\right) ,\frac{\partial u}{%
\partial t}\left( t\right) \right) ,\frac{\partial v}{\partial t}\left(
t\right) \right) \right] dt^{1}\wedge ...\wedge dt^{p}. 
\]
Moreover, 
\[
\left| \nabla _{x}L\left( t,u\left( t\right) ,\frac{\partial u}{\partial t}%
\left( t\right) \right) \right| \leq b_{0}\left( M\left| u\left( t\right)
\right| +\left| a\left( 0\right) \right| \right) \in L^{1}\left(
T_{0},R^{+}\right) 
\]
and 
\[
\left| \nabla _{y}L\left( t,u\left( t\right) ,\frac{\partial u}{\partial t}%
\left( t\right) \right) \right| \leq b_{0}\left( M\left| \frac{\partial u}{%
\partial t}\left( t\right) \right| +\left| a\left( 0\right) \right| \right)
\in L^{2}\left( T_{0},R^{+}\right) . 
\]
That is why

\[
\int_{T_{0}}\left[ \left( \nabla _{x}L\left( t,u\left( t\right) ,\frac{%
\partial u}{\partial t}\left( t\right) \right) ,v\left( t\right) \right)
\right. 
\]
\[
\left. +\left( \nabla _{y}L\left( t,u\left( t\right) ,\frac{\partial u}{%
\partial t}\left( t\right) \right) ,\frac{\partial v}{\partial t}\left(
t\right) \right) \right] dt^{1}\wedge ...\wedge dt^{p} 
\]

\[
\leq \int_{T_{0}}\left| \nabla _{x}L\left( t,u\left( t\right) ,\frac{%
\partial u}{\partial t}\left( t\right) \right) \right| \left| v\left(
t\right) \right| dt^{1}\wedge ...\wedge dt^{p} 
\]
\[
+\int_{T_{0}}\left| \nabla _{y}L\left( t,u\left( t\right) ,\frac{\partial u}{%
\partial t}\left( t\right) \right) \right| \left| \frac{\partial v}{\partial
t}\left( t\right) \right| dt^{1}\wedge ...\wedge dt^{p} 
\]
\[
\leq b_{0}\int_{T_{0}}\left( M\left| u\left( t\right) \right| +\left|
a\left( 0\right) \right| \right) \left| v\left( t\right) \right|
dt^{1}\wedge ...\wedge dt^{p} 
\]
\[
+b_{0}\int_{T_{0}}\left( M\left| \frac{\partial u}{\partial t}\left(
t\right) \right| +\left| a\left( 0\right) \right| \right) \left| \frac{%
\partial v}{\partial t}\left( t\right) \right| dt^{1}\wedge ...\wedge dt^{p} 
\]
By using the inequality Cauchy-Schwartz we find 
\[
\begin{array}{l}
\left| \displaystyle\frac{\partial \Psi }{\partial \lambda }\left( 0\right)
\right| \leq b_{0}\left( \displaystyle\int_{T_{0}}\left( M\left| u\left(
t\right) \right| +\left| a\left( 0\right) \right| \right) ^{2}dt^{1}\wedge
...\wedge dt^{p}\right) ^{\frac{1}{2}} \\ 
\noalign{\medskip}\cdot \left( \displaystyle\int_{T_{0}}\left| v\left(
t\right) \right| ^{2}dt^{1}\wedge ...\wedge dt^{p}\right) ^{\frac{1}{2}} \\ 
\noalign{\medskip}+b_{0}\left( \displaystyle\int_{T_{0}}\left( M\left| %
\displaystyle\frac{\partial u}{\partial t}\left( t\right) \right| +\left|
a\left( 0\right) \right| \right) ^{2}dt^{1}\wedge ...\wedge dt^{p}\right) ^{%
\frac{1}{2}}\left( \displaystyle\int_{T_{0}}\left| \displaystyle\frac{%
\partial v}{\partial t}\left( t\right) \right| ^{2}dt^{1}\wedge ...\wedge
dt^{p}\right) ^{\frac{1}{2}} \\ 
\noalign{\medskip}\leq C_{1}\left( \displaystyle\int_{T_{0}}\left| v\left(
t\right) \right| ^{2}dt^{1}\wedge ...\wedge dt^{p}\right) ^{\frac{1}{2}%
}+C_{2}\left( \displaystyle\int_{T_{0}}\left| \displaystyle\frac{\partial v}{%
\partial t}\left( t\right) \right| ^{2}dt^{1}\wedge ...\wedge dt^{p}\right)
^{\frac{1}{2}} \\ 
\noalign{\medskip}\leq \max \left\{ C_{1},C_{2}\right\} 2^{\frac{1}{2}%
}\left( \displaystyle\int_{T_{0}}\left( \left| v\left( t\right) \right|
^{2}+\left| \displaystyle\frac{\partial v}{\partial t}\left( t\right)
\right| ^{2}\right) dt^{1}\wedge ...\wedge dt^{p}\right) ^{\frac{1}{2}%
}=C\left\| v\right\| .
\end{array}
\]
By consequence, the action $\varphi $ has the derivative $\varphi ^{\prime
}\in \left( W_{T}^{1,2}\right) ^{\ast }$ given by $\left( 2\right) $. The
Krasnoselski theorem and the hypothesis $\left( 1\right) $ imply the fact
that the application $u\rightarrow \left( \nabla _{x}L\left( \cdot ,u,%
\displaystyle\frac{\partial u}{\partial t}\right) ,\nabla _{y}L\left( \cdot
,u,\displaystyle\frac{\partial u}{\partial t}\right) \right) ,$ from $%
W_{T}^{1,2}$ to $L^{1}\times L^{2},$ is continuous, so $\varphi ^{\prime }$
is continuous from $W_{T}^{1,2}$ to $\left( W_{T}^{1,2}\right) ^{\ast }$ and
the proof is complete.

\bigskip {\bf Theorem 3.} {\it If the action 
\[
\varphi \left( u\right) =\displaystyle\int_{T_{0}}\left[ \displaystyle\frac{1%
}{2}\left| \displaystyle\frac{\partial u}{\partial t}\right| ^{2}+F\left(
t,u\left( t\right) \right) \right] dt^{1}\wedge ...\wedge dt^{p} 
\]
is continuously differentiable on} ${\it H}_{{\it T}}^{{\it 1}}$ {\it and} $%
{\it u\in H}_{{\it T}}^{{\it 1}}${\it \ is a solution of the equation }${\it %
\varphi }^{{\it \prime }}{\it (u)=0}${\it \ (critical point), then the
function }${\it u}${\it \ has a weak Laplacian }${\it \bigtriangleup u}${\it %
\ (or the Jacobian matrix} ${\it \displaystyle\frac{\partial u}{\partial t}}$
{\it has a weak divergence) and } 
\[
{\it \bigtriangleup u=\nabla F(t,u(t))} 
\]
{\it a.e. on }${\it T}_{{\it 0}}${\it \ and } 
\[
{\it u\mid }_{{\it S}_{{\it i}}^{{\it -}}}{\it =u\mid }_{{\it S}_{{\it i}}^{%
{\it +}}}{\it ,\frac{\partial u}{\partial t}\mid }_{{\it S}_{{\it i}}^{{\it -%
}}}{\it =\frac{\partial u}{\partial t}\mid }_{{\it S}_{{\it i}}^{{\it +}}}%
{\it .} 
\]
{\bf Proof.} We build the function 
\[
\Phi :[-1,1]\rightarrow R, 
\]
\[
\Phi \left( \lambda \right) =\varphi \left( u+\lambda v\right) = 
\]
\[
\int_{T_{0}}\left[ \frac{1}{2}\left| \frac{\partial }{\partial t}\left(
u\left( t\right) +\lambda v\left( t\right) \right) \right| ^{2}+F\left(
t,u\left( t\right) +\lambda v\left( t\right) \right) \right] dt^{1}\wedge
...\wedge dt^{p}, 
\]
where $v\in C_{T}^{\infty }\left( T_{0},R^{n}\right) .$ The point $\lambda
=0 $ is a critical point of $\Phi $ if and only if the point $u$ is a
critical point of $\varphi $. Consequently 
\[
0=\left\langle \varphi ^{\prime }\left( u\right) ,v\right\rangle
=\int_{T_{0}}\left[ \delta ^{\alpha \beta }\delta _{ij}\frac{\partial u^{i}}{%
\partial t^{\alpha }}\frac{\partial v^{j}}{\partial t^{\beta }}+\delta
_{ij}\nabla ^{i}F\left( t,u\left( t\right) \right) v^{j}\left( t\right) %
\right] dt^{1}\wedge ...\wedge dt^{p}, 
\]
for all $v\in H_{T}^{1}$ and hence for all $v\in C_{T}^{\infty }.$ The
definition of the weak divergence, 
\[
\int_{T_{0}}\delta ^{\alpha \beta }\delta _{ij}\frac{\partial u^{i}}{%
\partial t^{\alpha }}\frac{\partial v^{j}}{\partial t^{\beta }}dt^{1}\wedge
...\wedge dt^{p} =-\int_{T_{0}}\delta ^{\alpha \beta }\delta _{ij}\frac{%
\partial ^{2}u^{i}}{\partial t^{\alpha }\partial t^{\beta }}%
v^{j}dt^{1}\wedge ...\wedge dt^{p}, 
\]
shows that the Jacobian matrix $\displaystyle\frac{\partial u}{\partial t}$
has weak divergence (the function $u$ has a weak Laplacian) and 
\[
\bigtriangleup u\left( t\right) =\nabla F\left( t,u\left( t\right) \right) 
\]
a.e. on $T_{0}$. Also, the existence of weak derivatives $\displaystyle\frac{%
\partial u}{\partial t}$ and weak divergence $\bigtriangleup u$ implies that 
\[
u\mid _{S_{i}^{-}}=u\mid _{S_{i}^{+}},\frac{\partial u}{\partial t}\mid
_{S_{i}^{-}}=\frac{\partial u}{\partial t}\mid _{S_{i}^{+}}. 
\]

\section{Poisson-gradient dynamical systems with periodical potential}

{\bf Theorem 4}. {\it If $F:T_{0}\times R^{n}\rightarrow R,$ $\left(
t,x\right) \rightarrow F\left( t,x\right) $ functions have the properties: }

{\it 1) $F\left( t,x\right) $ is measurable in $t$ for any $x\in R^{n}$ and
continuously differentiable in $x$ for any $t\in T_{0}$}${\it ,}${\it \ and
there exist the functions $a\in C^{1}\left( R^{+},R^{+}\right) $ with the
derivative $a^{\prime }$\ bounded from above and $b\in C\left(
T_{0},R^{+}\right) $ such that for any $t\in T_{0}$ and any $u\in R^{n}$ to
have $\left| F\left( t,u\right) \right| \leq a\left( \left| u\right| \right)
b\left( t\right) $ and $\left| \nabla _{u}F\left( t,u\right) \right| \leq
a\left( \left| u\right| \right) b\left( t\right) $}${\it ,}${\it \ }

{\it 2) $F\left( t,x\right) >0,$ for any $t\in T_{0}$ and any $x\in R^{n},$ }

{\it 3) For any $i\in \{1,...,n\}$ it exists $P_{i}\in R$ such that $F\left(
t,x+P_{i}e_{i}\right) =F\left( t,x\right) ,$ for any $t\in T_{0}$ and any $%
x\in R^{n},$ }

4) {\it The action} {\it $\varphi _{1}\left( u\right) =\displaystyle%
\int_{T_{0}}F\left( t,u\left( t\right) \right) dt^{1}\wedge ...\wedge dt^{p}$
is weakly lower semi-continous on }${\it H}_{{\it T}}^{{\it 1}}{\it .}$

\bigskip {\it If $\displaystyle\int_{T_{0}}F\left( t,u\right) dt^{1}\wedge
...\wedge dt^{p}\rightarrow \infty $ when $\left| u\right| \rightarrow
\infty $}${\it ,}${\it \ then, the problem $\Delta u\left( t\right) =\nabla
F\left( t,u\left( t\right) \right) $ with the boundary condition\ } {\it \ 
\[
u\mid _{S_{i}^{-}}=u\mid _{S_{i}^{+}},\displaystyle\frac{\partial u}{%
\partial t}\mid _{S_{i}^{-}}=\displaystyle\frac{\partial u}{\partial t}\mid
_{S_{i}^{+}} 
\]
} {\it has at least solution which minimizes the action 
\[
\varphi \left( u\right) =\displaystyle\int_{T_{0}}\left[ \displaystyle\frac{1%
}{2}\left| \displaystyle\frac{\partial u}{\partial t}\left( t\right) \right|
^{2}+F\left( t,u\left( t\right) \right) \right] dt^{1}\wedge ...\wedge
dt^{p} 
\]
in }${\it H}_{{\it T}}^{{\it 1}}${\it .}

{\it \noindent }{\bf Proof}. From Theorem 2, the action $\varphi $ is
continuously differentiable. From the periodicity and the continuity of $F,$
it results that exists the function $d\in $\ $L^{1}\left( T_{0},R\right) $
such that $F\left( t,x\right) \geq d\left( t\right) \geq 0$, for any $t\in
T_{0}$ and any $x\in R^{n}$. By consequence $\displaystyle%
\int_{T_{0}}F\left( t,u\left( t\right) \right) dt^{1}\wedge ...\wedge
dt^{p}=C_{1}\geq 0$. It results the inequality $\varphi \left( u\right) \geq %
\displaystyle\int_{T_{0}}\displaystyle\frac{1}{2}\left| \displaystyle\frac{%
\partial u}{\partial t}\right| ^{2}dt^{1}\wedge ...\wedge dt^{p}-C_{1}$ for
any $u\in H_{T}^{1}$. As result $\displaystyle\inf_{u\in H_{T}^{1}}\varphi
\left( u\right) <\infty $. Because $\varphi \left( u\right) +C_{1}\geq %
\displaystyle\int_{T_{0}}\displaystyle\frac{1}{2}\left| \displaystyle\frac{%
\partial u}{\partial t}\right| ^{2}dt^{1}\wedge ...\wedge dt^{p}$, for any $%
u\in H_{T}^{1}$, we have the same inequality and for any $u=u_{k}$ where $%
\left( u_{k}\right) $ it is a minimizing sequence for $\varphi $ in $%
H_{T}^{1}$. So, we obtain 
$$
\displaystyle\int_{T_{0}}\displaystyle\frac{1}{2}\left| \displaystyle\frac{%
\partial u_{k}}{\partial t}\right| ^{2}dt^{1}\wedge ...\wedge dt^{p}\leq
C_{2},\quad \hbox{for any}\quad k\in N.\eqno(3)
$$
\ We consider $u_{k}=\overline{u}_{k}+\widetilde{u}_{k}$, where $\overline{u}%
_{k}=\displaystyle\frac{1}{T^{1}...T^{p}}\int_{T_{0}}u_{k}\left( t\right)
dt^{1}\wedge ...\wedge dt^{p}$. From the relation $\left( 3\right) $ and the
Wirtinger inequality we have 
$$
\left\| \widetilde{u}_{k}\right\| \leq C_{3},\quad k\in N,\quad C_{3}\geq 0.%
\eqno(4)
$$
On the other side, from the periodicity of $F$ we find that $\varphi \left(
u+P_{i}e_{i}\right) =\varphi \left( u\right) $, for any $i\in \{1,...n\}$
and any $u\in H_{T}^{1}.$ If the sequence $\left( u_{k}\right) $ is a
minimizing one for the action $\varphi $, then the sequence $\left(
u_{k}^{\ast }\right) ,$ 
\[
u_{k}^{\ast }=\left( \overline{u}_{k}^{1}+\widetilde{u}%
_{k}^{1}+k_{1}P_{1},...,\overline{u}_{k}^{n}+\widetilde{u}%
_{k}^{n}+k_{n}P_{n}\right) 
\]
is also a minimizing sequence, for any $k_{1},...,k_{n}\in Z$. Obviously, we
may choose $k_{i}\in Z,i=1,...,n$ such that 
\[
0\leq \quad \overline{u}_{k}^{i}+k_{i}P_{i}\leq P_{i},i=1,...,n.
\]
From the relations $\left( 3\right) $ and $\left( 4\right) $ it results that
the sequence $\left( u_{k}^{\ast }\right) $ is bounded, so the action $%
\varphi $ has a minimizing bounded sequence. By eventually passing to a
subsequence, we may consider that the $\left( u_{k}^{\ast }\right) $
sequence is weakly convergent with the limit $u$.

The Hilbert space $H_{T}^{1}$ is reflexive. By consequence, the sequence $%
\left( u_{k}^{\ast }\right) $ (or one of his subsequence) is weakly
convergent in $H_{T}^{1}$ with the limit $u$. Because 
\[
\varphi _{2}\left( u\right) =\int_{T_{0}}\delta _{ij}\delta ^{\alpha \beta }%
\frac{\partial u^{i}}{\partial t^{\alpha }}\left( t\right) \frac{\partial
v^{j}}{\partial t^{\beta }}\left( t\right) dt^{1}\wedge ...\wedge dt^{p} 
\]
is convex it results that $\varphi _{2}$ is weakly lower semi-continuous, so 
\[
\varphi \left( u\right) =\varphi _{1}\left( u\right) +\varphi _{2}\left(
u\right) 
\]
is weakly lower semi-continuous and $\varphi \left( u\right) \leq \underline{%
\lim }\varphi \left( u_{k}\right) .$ This means that $u$ is minimum point of 
$\varphi .$ From the Theorem 3 this means that $u$ is solution of boundary
value problem 
\[
\bigtriangleup u\left( t\right) =\nabla F\left( t,u\left( t\right) \right) , 
\]

\[
u\mid _{S_{i}^{-}}=u\mid _{S_{i}^{+}},\frac{\partial u}{\partial t}\mid
_{S_{i}^{-}}=\frac{\partial u}{\partial t}\mid _{S_{i}^{+}}. 
\]

\bigskip \centerline{\bf References}\ \ 

[1] I. Duca, A-M. Teleman, C. Udri\c{s}te: {\it Poisson-Gradient Dynamical
Systems with Convex Potential}, Proceedings of the 3-rd International
Colloquium '' Mathematics in Engineering and Numerical Physics '', 7-9
October, 2004, Bucharest.

[2] Iulian Duca, Constantin Udri\c ste: {\it Some Inequalities satisfied by
Periodical Solutions of Multi-Time Hamilton Equations}, The 5-th Conference
of Balkan Society of Geometers, August 29-Sept 2, 2005, Mangalia, Romania.

[3] M. Forger, C. Paufler, H. Romer, {\it The Poisson bracket for Poisson
forms in multisymplectic field theory},  Reviews in Mathematical Physics,
15, 7 (2003), 705-743.

[4] J. Mawhin, M. Willem: {\it Critical Point Theory and Hamiltonian Systems}%
, Springer-Verlag, 1989.

[5] I. V. Kanatchikov: {\it Geometric (pre)quantization in the
polysymplectic approach to field theory}, arXiv: hep-th/0112263 v3, 3 Jun
2002, 1-12; Differential Geometry and Its Applications, Proc. Conf., Opava
(Czech Republic), August 27-31, 2001, Silesian University Opava, 2002.

[6] C. Paufler, H. Romer: {\it De Donder-Weyl equations and multisymplectic
geometry}, arXiv: math-ph/0506022 v2, 20 Jun 2005.

[7] N. Roman-Roy, {\it Multisymplectic Lagrangian and Hamiltonian formalism
of first -order classical field theories}, arXiv: math-ph/0107019 v1, 20 Jul
2001, vol XX (XXXX), No. X, 1-9.

[8] C.Udri\c{s}te: {\it Nonclassical Lagrangian Dynamics and Potential Maps}%
, Conference in Mathematics in Honour of Professor Radu Ro\c{s}ca on the
Occasion of his Ninetieth Birthday, Katholieke University Brussel,
Katholieke University Leuven, Belgium, Dec.11-16, 1999; \ 

http://xxx.lanl.gov/math.DS/0007060.

[9] C.Udri\c{s}te:{\it \ Solutions of DEs\ and PDEs\ as Potential Maps Using
First Order Lagrangians}, Centenial Vranceanu, Romanian Academy, University
of Bucharest, June 30-July 4, (2000); http://xxx.lanl.gov/math.DS/0007061;
Balkan Journal of Geometry and Its Applications 6, 1, 93-108, 2001.

[10] C. Udri\c{s}te, M. Postolache: {\it Atlas of Magnetic Geometric Dynamics%
}, Geometry Balkan Press, Bucharest, 2001.

[11] C. Udri\c{s}te: {\it From integral manifolds and metrics to potential
maps}, Atti del Academia Peloritana dei Pericolanti, Classe 1 di Science
Fis. Mat. e Nat., 81-82, A 01006 (2003-2004), 1-14.

[12] C. Udri\c{s}te, A-M. Teleman: {\it Hamilton Approaches of Fields Theory}%
, IJMMS, 57 (2004), 3045-3056; ICM Satelite Conference in Algebra and
Related Topics, University of Hong-Kong, 13-18.08.02.

[13] C. Udri\c{s}te, I. Duca: {\it Periodical Solutions of Multi-Time
Hamilton Equations,} Analele Universita\c{t}ii Bucure\c{s}ti, 55, 1 (2005),
179-188.

[14] C. Udri\c{s}te, I. Duca: \ {\it Poisson-gradient \ Dynamical \ Systems
with Bounded Non-linearity}, manuscript, 2005.

\bigskip

University Politehnica of Bucharest

Department of Mathematics

Splaiul Independentei 313

060042 Bucharest, Romania

email: udriste@mathem.pub.ro

\end{document}